\documentclass[12pt]{article}

\usepackage {amssymb}
\usepackage {amscd}
\usepackage {latexsym}
\usepackage {euscript}

\usepackage[T2A]{fontenc}
\usepackage[cp1251]{inputenc}
\usepackage[english, russian]{babel}
\usepackage[tbtags]{amsmath}
\usepackage{amsfonts,amssymb}

\usepackage{mathrsfs}
\usepackage{graphicx}

\newtheorem{theorem}{Theorem}
\newtheorem{definition}{Definition}
\newtheorem{proposition}{Proposition}

\newtheorem{remark}{Remark}

\addto{\captionsenglish}{}

\begin{document}
{\large
\begin{center} {\bfseries Multipliers in the Bessel potential spaces with positive smoothness indices: bilateral continuous embeddings and their exact character}
 \end{center}}

\bigskip
\centerline{ALEXEI\,A.~BELYAEV\footnote{This work is supported by Russian Science Foundation (RNF) under grant No 20-11-20261.}}

\bigskip

\medskip




\begin{quote}

{\large  Abstract.} We investigate the problem of establishing bilateral embeddings of the uniformly localized Bessel potential spaces $H^{\gamma}_{r, \: unif}(\mathbb{R}^n)$ into the multiplier spaces between Bessel potential spaces with positive smoothness indices. This problem is considered in the model situation when the natural norms of both of these Bessel potential spaces are generated by some inner product yet the description theorems for the corresponding multiplier space in terms of the spaces $H^{\gamma}_{r, \: unif}(\mathbb{R}^n)$ can not be established. The optimal character of the indices figuring in this embedding is also examined.

\end{quote}

\bigskip

\bigskip

\medskip


\bigskip

{\centerline {\Large 1. Introduction}}


\medskip

In this paper we consider the spaces of multipliers, acting in the scale of Bessel potential space $H^{\gamma}_r(\mathbb{R}^n)$, namely, we study the multiplier spaces $M[H^s_p(\mathbb{R}^n) \to H^t_q(\mathbb{R}^n)]$ in the situation when we impose the natural assumptions that guarantee the validity of the continuous embedding $H^s_p(\mathbb{R}^n) \underset{\to}{\subset} H^t_q(\mathbb{R}^n)$ and both smoothness indices $s$ and $t$ are positive and strictly less than $\frac{n}{2}$.

Our aim is twofold. Firstly, we want to establish bilateral continuous embeddings which involve the multiplier space $M[H^s_p(\mathbb{R}^n) \to H^t_q(\mathbb{R}^n)]$ and some spaces from the scale of the uniformly localized Bessel potential spaces $\{H^{\gamma}_{r, \: unif}(\mathbb{R}^n) \: | \; \gamma \in \mathbb{R}, \: r \in (1; \: +\infty)\}$. Secondly, we are interested in establishing the optimal character of the lower indices of the spaces $H^{\gamma}_{r, \: unif}(\mathbb{R}^n)$ featuring in these continuous embeddings.

For the situation $s \geqslant t > 0$ under the natural assumption $s > \frac{n}{p}\,$, guaranteeing the fact that $H^s_p(\mathbb{R}^n)$ can be can be endowed with the naturally defined multiplication which turns $H^s_p(\mathbb{R}^n)$ into an algebra, the studying of the sufficient conditions under which we can describe the multiplier space $M[H^s_p(\mathbb{R}^n) \to H^t_q(\mathbb{R}^n)]$ in terms of the uniformly localized Bessel potential spaces has quite a long history. It dates back to the seminal paper \cite{Str67} by R.\,S.~Stirchartz, where the uniformly localized Bessel potential spaces were initially introduced (though neither this term, nor notation $H^{\gamma}_{r, \: unif}(\mathbb{R}^n)$ were used in that paper). In that paper the following fact was established:

\medskip

{\bfseries Theorem.} \cite[Corollary 2.2]{Str67} {\itshape Let $p \in (1; \: +\infty)$ and let also $s \in \left(\frac{n}{p}; \; +\infty\right)$. Then the spaces $M[H^s_p(\mathbb{R}^n) \to H^s_p(\mathbb{R}^n)]$ and $H^s_{p, \: unif}(\mathbb{R}^n)$ coincide with their natural norms being equivalent.}

\bigskip

The situation when the smoothness indices can be different was treated much later (see, e.g., \cite[Chapter 3, Theorem 3.2.5]{MShbook}) and it turned out that if $p \in (1; \: +\infty), \: s > t \geqslant 0$ and $s > \frac{n}{p}$, then $M[H^s_p(\mathbb{R}^n) \to H^t_p(\mathbb{R}^n)] = H^t_{p, \: unif}(\mathbb{R}^n)$ and the norms of these spaces are equivalent.

The most general situation of the multiplier space $M[H^s_p(\mathbb{R}^n) \to H^t_q(\mathbb{R}^n)]$, when lower indices $p$ and $q$ may also differ, was examined by the author and A.\,A.~Shkalikov in \cite{BelShk2017}, where under some natural assumptions on indices $s, \: t, \: p$ and $q$ the following result was proved:

\medskip

{\bfseries Theorem.} \cite[Theorem 1]{BelShk2017} {\itshape Let $p \in (1; \: +\infty), \; q \in (1; \: +\infty), \; p \leqslant q, \; s \in \left(\frac{n}{p}; \: +\infty\right), \; t \in [0; \: +\infty)$ and $s - \frac{n}{p} \geqslant t - \frac{n}{q}$. Then the spaces $M[H^s_p(\mathbb{R}^n) \to H^t_q(\mathbb{R}^n)]$ and $H^t_{q, \: unif}(\mathbb{R}^n)$ coincide with their natural norms being equivalent.}

\bigskip

It is worthy to note that the similar problem of finding a description of the multiplier space $M[A^s_{p_1, \: q_1}(\mathbb{R}^n) \to A^t_{p_2, \: q_2}(\mathbb{R}^n)]$ in terms of the uniformly localized spaces $A^{\gamma}_{p, \: q, \: unif}(\mathbb{R}^n)$ in the situation when $A$ stands either for Lizorkin--Triebel spaces (of which Bessel potential spaces $H^s_p(\mathbb{R}^n)$ present a partial case as for any $s \in \mathbb{R}$ and $p \in [1; \: +\infty)$ we have $H^s_p(\mathbb{R}^n) = F^s_{p, \: 2}(\mathbb{R}^n)$) or Besov spaces (or some anisotropic versions of the Besov--Lizorkin--Triebel type spaces) and the Strichartz--type condition $s > \frac{n}{p}$ holds true, was extensively studied and numerous results, analogous to Strichartz's theorem, were obtained (see \cite{Bour, Fr, NguyenSickel} and also monograph \cite{RSbook} and the references therein).

The case of smoothness indices of different signs, which is very important for its many fruitful applications to the spectral theory of differential operators, including the study of the general type singular perturbations of the strongly elliptic operators and the Navier--Stokes problem (see, e.g., monographs \cite{AGHH2ndEd} and \cite{Lemarie-Rieusset_NS_21}, and also \cite{KMPaper, LRG, GerPaper, NZSh2006}), was extensively studied in the series of papers by A.\,A.~Shkalikov, M.\,I.~Neiman--Zade, J.\,G.~Bak and the author (see \cite{BSh, Bel2, BelShk2018, NZSh2006}). The most general result in this direction on the description of the multiplier space $M[H^s_p(\mathbb{R}^n) \to H^{-t}_q(\mathbb{R}^n)]$ for nonnegative $s$ and $t$ was obtained in \cite{BelShk2018}. It reads as follows:

\medskip

{\bfseries Theorem.} \cite[Theorem 1]{BelShk2018} {\itshape Let $p \in (1; \: +\infty), \; q \in (1; \: +\infty), \; p \leqslant q$ and either $s \geqslant t \geqslant 0$ and $s > \frac{n}{p}$, or $t \geqslant s \geqslant 0$ and $t > \frac{n}{q'}$. Then the spaces $M[H^s_p(\mathbb{R}^n) \to H^{-t}_q(\mathbb{R}^n)]$ and $H^{-t}_{q, \: unif}(\mathbb{R}^n) \cap H^{-s}_{p', \; unif}(\mathbb{R}^n)$ coincide with their natural norms being equivalent.}

\bigskip

So it turns out that in the case when some generalizations of the condition $s > \frac{n}{p}$ do hold, under natural additional assumptions we can obtain the constructive characterization of the multiplier spaces $M[H^s_p(\mathbb{R}^n) \to H^t_q(\mathbb{R}^n)]$, no matter if both smoothness indices $s$ and $t$ are nonnegative or $s$ is nonnegative and $t$ is nonpositive. This naturally leads us to the following question: what is the relation between the multiplier spaces $M[H^s_p(\mathbb{R}^n) \to H^t_q(\mathbb{R}^n)]$ and the uniformly localized Bessel potential spaces in the situation when Strichartz--type conditions do not hold?

The model case $p = q = 2$ (in which the norms $\| \cdot \|_{H^s_p(\mathbb{R}^n)}$ and $\| \cdot \|_{H^t_q(\mathbb{R}^n)}$ can be induced by the inner products) for the situation when smoothness indices are of different sign was treated in \cite{NZSh2006} and \cite{BelShkArXiv2019}. Namely, the following results were established.

\medskip

{\bfseries Theorem.} (An immediate corollary of \cite[Lemma 3, Lemma 6]{NZSh2006}) {\itshape Let $s \in [0; \: +\infty]$, $ t \in [0; \: +\infty)$ and $\max(s; \: t) \in \left(0; \: \frac{n}{2}\right)$.

Then the following continuous embeddings hold true:
$$
H^{-\min(s; \: t)}_{\frac{n}{\max(s; \: t)}, \: unif}(\mathbb{R}^n) \underset{\to}{\subset} M[H^s_2(\mathbb{R}^n) \to H^{-t}_2(\mathbb{R}^n)]
$$
\centerline{and}
$$
M[H^s_2(\mathbb{R}^n) \to H^{-t}_2(\mathbb{R}^n)] \underset{\to}{\subset} H^{-\min(s; \: t)}_{2, \: unif}(\mathbb{R}^n).
$$}

\bigskip

In the paper \cite{BelShkArXiv2019} the exact character of the left embedding in this theorem was established by considering a special class of regular distributions. Namely, the following result about optimality of the lower index $\frac{n}{\max (s, t)}$ in the continuous embedding
$$
H^{-\min(s; \: t)}_{\frac{n}{\max(s; \: t)}, \: unif}(\mathbb{R}^n) \underset{\to}{\subset} M[H^s_2(\mathbb{R}^n) \to H^{-t}_2(\mathbb{R}^n)]
$$
was proved.

\medskip

{\bfseries Theorem.} \cite[Theorem 3]{BelShkArXiv2019} {\itshape Let $s \in [0; \: +\infty), \: t \in [0; \: +\infty)$ and $\max (s, t) \in \left(0; \: \frac{n}{2}\right)$. Then for arbitrary $\varepsilon \in \left(0; \: \frac{n}{\max (s, t)} - 2\right)$ there exists a distribution $u_{\varepsilon} \in \mathcal{S}^{'}(\mathbb{R}^n)$, such that $u_{\varepsilon}$ belongs to $H^{-\min (s, t)}_{\frac{n}{\max (s, t)} - \varepsilon, \: unif}(\mathbb{R}^n)$ but it is not a multiplier from $H^s_2(\mathbb{R}^n)$ to $H^{-t}_2(\mathbb{R}^n)$, i. e.
$$
u_{\varepsilon} \in H^{-\min (s, t)}_{\frac{n}{\max (s, t)} - \varepsilon, \: unif}(\mathbb{R}^n) \setminus \: M[H^s_2(\mathbb{R}^n) \to H^{-t}_2(\mathbb{R}^n)].
$$}

\bigskip

The main aim of this paper is to obtain the analogues of these two theorems for the multiplier space $M[H^s_p(\mathbb{R}^n) \to H^t_q(\mathbb{R}^n)]$ in the case when both smoothness indices $s$ and $t$ are positive yet the condition $s > \frac{n}{p}$ does not hold true.

Namely, we consider the situation when $s \in (0; \: +\infty), \: t \in [0; \: +\infty), s > t$ and $s < \frac{n}{2}$.

Since it can be easily seen that we have the continuous embedding $H^t_{2, \: unif}(\mathbb{R}^n) \underset{\to}{\subset} H^{-s}_{2, \: unif}(\mathbb{R}^n)$, the embedding $M[H^s_2(\mathbb{R}^n) \to H^t_2(\mathbb{R}^n)] \underset{\to}{\subset} H^t_{2, \: unif}(\mathbb{R}^n)$ is a trivial corollary of the general fact (see \cite[Proposition 2]{BelShk2017}) that for any $s \in \mathbb{R}, \: t \in \mathbb{R}, \: p \in (1; \: +\infty)$ and $q \in (1; \; +\infty)$ we have the continuous embedding
$$
M[H^s_p(\mathbb{R}^n) \to H^t_q(\mathbb{R}^n)] \underset{\to}{\subset} H^t_{q, \: unif}(\mathbb{R}^n) \cap H^{-s}_{p', \: unif}(\mathbb{R}^n).
$$

Continuous embedding $H^t_{\frac{n}{s}, \: unif}(\mathbb{R}^n) \underset{\to}{\subset} M[H^s_2(\mathbb{R}^n) \to H^t_2(\mathbb{R}^n)]$ is proved via the criterion of the continuous embedding
$$
H^{\gamma}_{r, \: unif}(\mathbb{R}^n) \underset{\to}{\subset} M[H^s_p(\mathbb{R}^n) \to H^t_q(\mathbb{R}^n)]
$$
in terms of the validity of the multiplicative functional estimate, established in \cite{BelShk2017}, and the general multiplicative estimates for the Lizorkin--Triebel spaces from \cite{RSbook}.

By implementing this approach, we obtain the following result.

\medskip

{\bfseries Theorem 4.} {\itshape Let $s \in (0; \: +\infty), \: t \in [0; \: +\infty), s > t$ and $s < \frac{n}{2}$.

Then the continuous embeddings
$$
H^t_{\frac{n}{s}, \: unif}(\mathbb{R}^n) \underset{\to}{\subset} M[H^s_2(\mathbb{R}^n) \to H^t_2(\mathbb{R}^n)] \quad \mbox{and} \quad M[H^s_2(\mathbb{R}^n) \to H^t_2(\mathbb{R}^n)] \underset{\to}{\subset} H^t_{2, \: unif}(\mathbb{R}^n)
$$
hold true.}

\bigskip

The exact character of the left embedding $H^t_{\frac{n}{s}, \: unif}(\mathbb{R}^n) \underset{\to}{\subset} M[H^s_2(\mathbb{R}^n) \to H^t_2(\mathbb{R}^n)]$ is proved by considering the special class of regular distributions $\eta \cdot \mathbf{f}_{\alpha}$, where $\eta$ is a smooth function from $\mathcal{D}(\mathbb{R}^n)$, satisfying some additional assumptions, and the function $f_{\alpha} \colon \mathbb{R}^n \to \mathbb{R}$ is defined as follows:
$$
\forall \: x \in \mathbb{R}^n \; \: f_{\alpha}(x) = \begin{cases}
\frac{1}{\left|x\right|^{\alpha}}, \; \: \mbox{if} \: x \neq \mathbf{0};\\
0, \; \: \mbox{if} \: x = \mathbf{0}.
\end{cases}
$$

These considerations allow us to prove the following result.

\medskip

{\bfseries Theorem 5.} {\itshape Let $s \in \left(0; \: \frac{n}{2}\right), \: t \in \left(0; \: \frac{n}{2}\right)$ and $s > t$. Then for any number $\varepsilon \in \left(0; \; \frac{n}{s} - 2\right)$ there exists a distribution
$$
u_{\varepsilon} \in H^t_{\frac{n}{s} - \varepsilon, \: unif}(\mathbb{R}^n) \setminus M[H^s_2(\mathbb{R}^n) \to H^t_2(\mathbb{R}^n)].
$$}

\bigskip

Summarizing the above considerations, we can make a conclusion that for both cases, when smoothness indices $s$ and $t$ are positive and when $s$ is positive and $t$ is negative, under natural additional assumptions (without which no uniformly localized Bessel potential space is continuously embedded into the multiplier space $M[H^s_2(\mathbb{R}^n) \to H^t_2(\mathbb{R}^n)]$, we can establish the bilateral continuous embeddings of the type
$$
H^{\gamma}_{r_1; \: unif}(\mathbb{R}^n) \subset M[H^s_2(\mathbb{R}^n) \to H^t_2(\mathbb{R}^n)] \subset H^{\gamma}_{r_2; \: unif}(\mathbb{R}^n)
$$
with the left of these embeddings being exact in the sense that its analogue is no longer valid if we decrease (no matter how little) the index $r_1$.

The major importance of this result about the exactness of the left--hand side embedding lies in the fact that it implies that, once Strichartz--type condition is no longer valid, we can not obtain the analogues of the description theorems for the multiplier space $M[H^s_p(\mathbb{R}^n) \to H^t_q(\mathbb{R}^n)]$ in terms of their coincidence with some uniformly localized Bessel potential space (or the intersection of those) even in the most basic situation when $p = q = 2$, no matter if the smoothness indices of these Bessel potential spaces $H^s_p(\mathbb{R}^n)$ and $H^t_q(\mathbb{R}^n)$ are of the same sign or not.

\bigskip

\bigskip

\bigskip

\bigskip

{\centerline {\Large 2. Basic definitions and classical results: spaces of test functions}}

\medskip

{\centerline {\Large and distributions, scales of the Bessel potential spaces}}

\medskip

{\centerline {\Large and of the uniformly localized Bessel potential spaces.}}

\bigskip

\bigskip




In what follows, we shall say that the normed space $(X, \: \| \cdot \|_X)$ is continuously embedded into the normed space $(Y, \: \| \cdot \|_Y)$ and denote this by $X \underset{\to}{\subset} Y$, if $X \subset Y$ and there exists a constant $C \in (0; \: +\infty)$, such that
$$
\forall \: u \in X \; \; \| u \|_Y \leqslant C \cdot \| u \|_X.
$$

By $\mathbb{Z}_+$ we shall denote the set of all nonnegative integer numbers and for multi-index $\beta = \left(\beta_k \: | \; k \in \overline{1, \: n}\right) \in \mathbb{Z}^n$ we shall use the notation $|\beta|_1 = \sum\limits_{k = 1}^n |\beta_k|$. If a function $f \colon \mathbb{R}^n \to \mathbb{C}$ is infinitely differentiable on $\mathbb{R}^n$ (in a classical sense), then for arbitrary multi-index $\alpha = \left(\alpha_k \: | \; k \in \overline{1, \: n}\right) \in \left(\mathbb{Z}_+\right)^n$ by $D^{\alpha}(f)$ we shall denote the action of the differential operator $D^{\alpha}$ on $f$, where an order of $D^{\alpha}$ equals $|\alpha|_1$ and for any $k \in \overline{1, \: n}$ order of the differentiation by $k^{\mbox{th}}$ variable equals $\alpha_k$.

By $\mu_n$ we shall denote the $n$--dimensional classical Lebesgue measure, and for $p \in (1; \; +\infty)$ by $p'$ we shall denote its Lebesgue dual number $p' \in (1; \: +\infty)$, defined by equality
$
\frac{1}{p} + \frac{1}{p'} = 1.
$

Let $p \in [1; \: +\infty)$. Then $\mathcal{L}_p(\mathbb{R}^n)$ is a complex linear space (endowed with the pointwise linear operations) of all Lebesgue measurable functions $f \colon \mathbb{R}^n \to \mathbb{C}$, such that the function $\left|f\right|^p$ is integrable with respect to the classical Lebesgue measure $\mu_n$. This linear space is endowed with a natural seinorm $\| \cdot \|_{\mathcal{L}_p} \colon \mathcal{L}_p(\mathbb{R}^n) \to \mathbb{R}$, defined by
$$
\forall \: f \in \mathcal{L}_p(\mathbb{R}^n) \; \; \; \| f \|_{\mathcal{L}_p} = \left(\;\int\limits_{\mathbb{R}^n}  \left|f\right|^p \: d\mu_n \right)^{\frac{1}{p}}.
$$
By $L_p(\mathbb{R}^n)$ we denote the linear space which is the result of taking the quotient space of $\mathcal{L}_p(\mathbb{R}^n)$ by the kernel of its seminorm $\| \cdot \|_{\mathcal{L}_p}$, i.e. the set $L_p(\mathbb{R}^n)$ consists of the classes of equivalent functions belonging to $\mathcal{L}_p(\mathbb{R}^n)$ with two functions $f_1 \in \mathcal{L}_p(\mathbb{R}^n)$ and $f_2 \in \mathcal{L}_p(\mathbb{R}^n)$ belonging to the same equivalence class if and only if $\mu_n\left(\{x \in \mathbb{R}^n \: | \; f_1(x) \neq f_2(x)\}\right) = 0$. On the linear space $L_p(\mathbb{R}^n)$ the norm $\| \cdot \|_{L_p} \colon L_p(\mathbb{R}^n) \to \mathbb{R}$ can be naturally introduced with its value on the arbitrary equivalence class from $L_p(\mathbb{R}^n)$ being equal to the value of the seminorm $\| \cdot \|_{\mathcal{L}_p}$ on arbitrary function from this equivalence class.

By $\mathcal{D}(\mathbb{R}^n)$ and $\mathcal{S}(\mathbb{R}^n)$ we denote the standard linear spaces of the test functions, consisting respectively of all the infinitely differentiable functions with compact support and all the infinitely diferentiable functions which are rapidly decreasing in the Schwartz sense. On the linear space $\mathcal{S}(\mathbb{R}^n)$ topology $\tau_{\mathcal{S}}$, which endows $\mathcal{S}(\mathbb{R}^n)$ with a structure of a linear topological space, is defined via the countable family of seminorms $\left\{\| \cdot \|_{k, \: m} \: | \; k \in \mathbb{Z}_+, \; m \in \mathbb{Z}_+ \right\}$, where for arbitrary numbers $k \in \mathbb{Z}_+$ and $m \in \mathbb{Z}_+$ we define the seminorm $\| \cdot \|_{k, \: m}$ by letting
$$
\forall \: f \in \mathcal{S}(\mathbb{R}^n) \; \; \; \| f \|_{k, \: m} = \max\limits_{\alpha \in \left(\mathbb{Z}_+\right)^n \colon |\alpha|_1 \leqslant m} \left(\sup\limits_{x \in \mathbb{R}^n} \left(\left(1 + \left|x\right|^2\right)^{\frac{k}{2}} \cdot \left|(D^{\alpha}(f))(x)\right| \right)\right).
$$
The space $\mathcal{D}(\mathbb{R}^n)$ is also a linear topological space with respect to the naturally introduced topology $\tau_{\mathcal{D}}$, such that the sequence of functions $\left(f_k \in \mathcal{D}(\mathbb{R}^n) \: | \; k \in \mathbb{N}\right)$ converges to a function $f_0 \in \mathcal{D}(\mathbb{R}^n)$ with respect to the topology $\tau_{\mathcal{D}}$ if and only if there exists a compact $K \subset \mathbb{R}^n$, such that the supports of all functions $f_k, \: k \in \mathbb{Z}_+$, are subsets of $K$ and for any multi-index $\alpha \in \left(\mathbb{Z}_+\right)^n$ the sequence of functions $\left(D^{\alpha}(f_k) \: | \; k \in \mathbb{N}\right)$ converges to the funcion $D^{\alpha}(f_0)$ uniformly on $K$. It is well--known (see, e.g., \cite[pp. 333 -- 334]{Khelem}) that this topology $\tau_{\mathcal{D}}$ is not countably (semi)normable but nevertheless it can be induced by some (uncountable) family of norms.

By $\mathcal{D}^{'}(\mathbb{R}^n)$ and $\mathcal{S}^{'}(\mathbb{R}^n)$ we shall denote the spaces of distributions, dual to the spaces $\mathcal{D}(\mathbb{R}^n)$ and $\mathcal{S}(\mathbb{R}^n)$ respectively, which consist of all sequentially continuous (with respect to the topologies $\tau_{\mathcal{D}}$ and $\tau_{\mathcal{S}}$ respectively) linear functionals, defined on the corresponding spaces of test functions. Linear spaces $\mathcal{D}^{'}(\mathbb{R}^n)$ and $\mathcal{S}^{'}(\mathbb{R}^n)$ are endowed with the topologies $\tau_{\mathcal{D}^{'}}$ and $\tau_{\mathcal{S}^{'}}$ of $*$--weak convergence, induced by the families of seminorms $\{\left.\| \cdot \|_f \: \right| f \in \mathcal{D}(\mathbb{R}^n)\}$ and $\{\left.\| \cdot \|_g \: \right| g \in \mathcal{S}(\mathbb{R}^n)\}$ respectively, where
$$
\forall \: f \in \mathcal{D}(\mathbb{R}^n) \; \; \mbox{we let} \; \; \forall \: u \in \mathcal{D}'(\mathbb{R}^n) \; \; \| u \|_f = |u(f)|
$$
\centerline{and}
$$
\forall \: g \in \mathcal{S}(\mathbb{R}^n) \; \; \mbox{we let} \; \; \forall \: v \in \mathcal{S}'(\mathbb{R}^n) \; \; \| v \|_g = |v(g)|.
$$
Also for an arbitrary number $p \in [1; \: +\infty)$ by $\mathcal{L}_{p, \: loc}(\mathbb{R}^n)$ we shall denote the linear space of all functions $f \colon \mathbb{R}^n \to \mathbb{C}$, such that for arbitrary function $\varphi \in \mathcal{D}(\mathbb{R}^n)$ we have $\varphi \cdot f \in \mathcal{L}_p(\mathbb{R}^n)$.

Let us also remind that for an arbitrary function $f \in \mathcal{L}_1(\mathbb{R}^n)$ its Fourier transform is a function $\mathcal{F}(f) \colon \mathbb{R}^n \to \mathbb{C}$, defined by
$$
\forall \: x \in \mathbb{R}^n \; \: (\mathcal{F}(f))(x) = \frac{1}{\left(2 \cdot \pi\right)^{\frac{n}{2}}} \cdot \int\limits_{\mathbb{R}^n} e^{-i \cdot < x \, ; \: y >} \cdot f(y) \: d\mu_n(y),
$$
and its inverse Fourier transform is the function $\mathcal{F}^{-1}(f) \colon \mathbb{R}^n \to \mathbb{C}$, defined by
$$
\forall \: x \in \mathbb{R}^n \; \: (\mathcal{F}^{-1}(f))(x) = \frac{1}{\left(2 \cdot \pi\right)^{\frac{n}{2}}} \cdot \int\limits_{\mathbb{R}^n} e^{i \cdot < x \, ; \: y >} \cdot f(y) \: d\mu_n(y).
$$
As it is well--known (see, e.g., \cite[Theorem 2.2.14 and Corollary 2.2.15]{Grafakos_Cl}), the Fourier transform is a homeomorphism (with repect to the topology $\tau_{\mathcal{S}}$) of a topological linear space $\mathcal{S}(\mathbb{R}^n)$ onto itself and inverse of the mapping $\mathcal{F} \colon \mathcal{S}(\mathbb{R}^n) \to \mathcal{S}(\mathbb{R}^n)$ is the mapping $\mathcal{F}^{-1} \colon \mathcal{S}(\mathbb{R}^n) \to \mathcal{S}(\mathbb{R}^n)$. This fact allows us to define both the Fourier transform $\mathcal{F}$ and the inverse Fourier tranform $\mathcal{F}^{-1}$ on the dual Schwartz space $\mathcal{S}'(\mathbb{R}^n)$ as follows:

for an arbitrary distribution $u \in \mathcal{S}'(\mathbb{R}^n)$ the distribution $\mathcal{F}(u) \in \mathcal{S}'(\mathbb{R}^n)$ is defined by
$$
\forall \: \varphi \in \mathcal{S}(\mathbb{R}^n) \; \; \: (\mathcal{F}(u))(\varphi) = u(\mathcal{F}(\varphi)),
$$
and the distribution $\mathcal{F}^{-1}(u) \in \mathcal{S}'(\mathbb{R}^n)$ is defined by
$$
\forall \: \varphi \in \mathcal{S}(\mathbb{R}^n) \; \; \: (\mathcal{F}^{-1}(u))(\varphi) = u(\mathcal{F}^{-1}(\varphi)).
$$

It immediately follows that the mappings $\mathcal{F} \colon \mathcal{S}'(\mathbb{R}^n) \to \mathcal{S}'(\mathbb{R}^n)$ and $\mathcal{F}^{-1} \colon \mathcal{S}'(\mathbb{R}^n) \to \mathcal{S}'(\mathbb{R}^n)$ are also linear homeomorphisms with respect to the topology $\tau_{\mathcal{S}'}$ and these mappings are inverses of each other.

Let us introduce the following notations for an open ball, a closed ball and a spherical shell which will be used throughout this paper. Let $r \in (0; \: +\infty), \: r_1 \in (0; \: +\infty), \: r_2 \in (r_1; \: +\infty)$ and $x_0 \in \mathbb{R}^n$. Then
$$
B_r(x_0) \stackrel{def}{=} \left\{ \left. x \in \mathbb{R}^n \: \right| \; \: |x - x_0| < r \right\}, \quad \overline{B_r}(x_0) \stackrel{def}{=} \left\{ \left. x \in \mathbb{R}^n \: \right| \; \: |x - x_0| \leqslant r \right\},
$$
$$
B_{r_1, \: r_2}(x_0) \stackrel{def}{=} \left\{ x \in \mathbb{R}^n \: \left| \; \: |x - x_0| > r_1 \; \: \mbox{and} \; \: |x - x_0| < r_2 \right. \right\}.
$$

Let $f \in \mathcal{L}_{1, \: loc}(\mathbb{R}^n)$ be a function such that it has at most polynomial growth outside of some ball, i.e. there exists such number $r_0 \in (0; \; +\infty)$, that for some numbers $N \in \mathbb{N}$ and $C \in [0; \: +\infty)$ the following inequality holds true:
$$
\forall \; x \in \mathbb{R}^n \setminus B_{r_0}(\mathbf{0}) \; \; \; \: |f(x)| \leqslant C \cdot \left(1 + \left|x\right|^2\right)^N.
$$
It is easy to see that for an arbitrary function $\varphi \in \mathcal{S}(\mathbb{R}^n)$ we have $\varphi \cdot f \in \mathcal{L}_1(\mathbb{R}^n)$. Indeed, this function is integrable on $\overline{B_{r_0}}(\mathbf{0})$ because $f \in \mathcal{L}_{1, \: loc}(\mathbb{R}^n)$ and, hence, function $f$ is integrable on the closed ball $\overline{B_{r_0}}(\mathbf{0})$, while function $\varphi$, which is infinitely differentiable on $\mathbb{R}^n$, is bounded on $\overline{B_{r_0}}(\mathbf{0})$. On the other hand, the integrability of $\varphi \cdot f$ outside the ball $\overline{B_{r_0}}(\mathbf{0})$ follows immediately from the integrability of the function
$$
g \colon x \longmapsto \frac{C \cdot \, \| \varphi \|_{2 \cdot (N + n), \: 0}}{\left(1 + \left|x\right|^2\right)^n}
$$
on $\mathbb{R}^n \setminus \overline{B_{r_0}}(\mathbf{0})$ and validity of the estimate
$$
\forall \; x \in \mathbb{R}^n \setminus \overline{B_{r_0}}(\mathbf{0}) \; \; \; \left|f(x) \cdot \varphi(x)\right| \leqslant C \cdot \left(1 + \left|x\right|^2\right)^N \cdot |\varphi(x)| =
$$
$$
= \frac{C}{\left(1 + \left|x\right|^2\right)^n} \cdot \left(1 + \left|x\right|^2\right)^{N + n} \cdot |\varphi(x)| \leqslant \frac{C \cdot \, \| \varphi \|_{2 \cdot (N + n), \: 0}}{\left(1 + \left|x\right|^2\right)^n} = g(x).
$$
Since for any function $\varphi \in \mathcal{S}(\mathbb{R}^n)$ we have
$$
\left| \; \int\limits_{\mathbb{R}^n} f(x) \cdot \varphi(x) \: d\mu_n(x) \right| \leqslant \int\limits_{B_{r_0}(\mathbf{0})} |f(x)| \cdot | \varphi(x)| \: d\mu_n(x) \; + \int\limits_{\mathbb{R}^n \setminus B_{r_0}(\mathbf{0})} |f(x)| \cdot |\varphi(x)| \: d\mu_n(x) \leqslant
$$
$$
\leqslant K_1 \cdot \|\varphi \|_{0, \: 0} + K_2 \cdot \| \varphi \|_{2 \cdot (N + n), \: 0} \leqslant (K_1 + K_2) \cdot \| \varphi \|_{2 \cdot (N + n), \: 0},
$$
where
$$
K_1 = \int\limits_{B_{r_0}(\mathbf{0})} |f(x)| \: d\mu_n(x), \quad K_2 = C \cdot \int\limits_{\mathbb{R}^n \setminus B_{r_0}(\mathbf{0})} \frac{1}{\left(1 + \left|x\right|^2\right)^n} \: d\mu_n(x),
$$
the regular functional $\mathbf{f}$, generated by the function $f$ and acting from $\mathcal{S}(\mathbb{R}^n)$ to $\mathbb{C}$, defined by
\begin{equation}\label{regular_eq}
\forall \: \varphi \in \mathcal{S}(\mathbb{R}^n) \; \: \; \; \; \: \mathbf{f}(\varphi) = \int\limits_{\mathbb{R}^n} f \cdot \varphi \: d\mu_n,
\end{equation}
is well--defined. From the linearity of this functional $\mathbf{f}$ and its continuity with respect to $\tau_{\mathcal{S}}$ it follows that $\mathbf{f}$ belongs to the space $\mathcal{S}'(\mathbb{R}^n)$. The restriction of this functional on $\mathcal{D}(\mathbb{R}^n)$, which will be also denoted by $\mathbf{f}$ where it does not lead to confusion, also turns out to be sequentially continuous with respect to $\tau_{\mathcal{D}}$ and, therefore, it belongs to the space $\mathcal{D}'(\mathbb{R}^n)$.

Let $p \in [1; \: +\infty)$. Since for any two functions $f \in \mathcal{L}_p(\mathbb{R}^n)$ and $g \in \mathcal{L}_p(\mathbb{R}^n)$, differing on the null set (with respect to the classical Lebesgue measure $\mu_n$), distributions $\mathbf{f}$ and $\mathbf{g}$ coincide as elements of $\mathcal{S}'(\mathbb{R}^n)$ (and, hence, as elements of $\mathcal{D}'(\mathbb{R}^n)$), mapping a class of equivalent functions, generated by the function $f \in \mathcal{L}_p(\mathbb{R}^n)$, to the distribution $\mathbf{f}$ actually implements an embedding of $L_p(\mathbb{R}^n)$ into $\mathcal{S}'(\mathbb{R}^n)$ and into $\mathcal{D}'(\mathbb{R}^n)$.

Let $\psi \colon \mathbb{R}^n \to \mathbb{C}$ be an arbitrary function. Then 
for arbitrary $y \in \mathbb{R}^n$ we shall denote by $\psi_{(y)}$ its shift by $y$, i.e. the function $\psi_{(y)} \colon \mathbb{R}^n \to \mathbb{C}$, defined by
$$
\forall \: x \in \mathbb{R}^n \; \: \psi_{(y)}(x) = \psi(x - y).
$$

Let $g \in \mathcal{S}(\mathbb{R}^n)$ and $g_1 \in \mathcal{S}(\mathbb{R}^n)$ be such functions that the equality $\mathcal{F}(g) = g_1$ holds true. From the basic fact of harmonic analysis (see, e.g., \cite[Theorem 2.2.14 (1)]{Grafakos_Cl}), which states that for any two functions $\varphi_1 \in \mathcal{S}(\mathbb{R}^n)$ and $\varphi_2 \in \mathcal{S}(\mathbb{R}^n)$ we have
$$
\int\limits_{\mathbb{R}^n} \mathcal{F}(\varphi_1) \cdot \varphi_2 \: d\mu_n = \int\limits_{\mathbb{R}^n} \varphi_1 \cdot \mathcal{F}(\varphi_2) \: d\mu_n,
$$
the following chain of equalities immediately follows for an arbitrary function $\psi \in \mathcal{S}(\mathbb{R}^n)$:
$$
\left(\mathcal{F}\left(\mathbf{g}\right)\right)(\psi) = \mathbf{g}\left(\mathcal{F}(\psi)\right) = \int\limits_{\mathbb{R}^n} g \cdot \mathcal{F}(\psi) \: d\mu_n = \int\limits_{\mathbb{R}^n} \mathcal{F}(g) \cdot \psi \: d\mu_n = \int\limits_{\mathbb{R}^n} g_1 \cdot \psi \: d\mu_n = \mathbf{g}_1(\psi).
$$
This means that the notions of the Fourier transform (and, hence, the notions of the inverse Fourier transform), defined on $\mathcal{S}(\mathbb{R}^n)$ and on $\mathcal{S}'(\mathbb{R}^n)$, agree with each other.

Let us now introduce the scale of the Bessel potential spaces and the scale of the uniformly localized Bessel potential spaces.

\begin{definition}\label{def_H^0_p}
Let $p \in [1; \: +\infty)$. Let also $v \in \mathcal{S}^{'}(\mathbb{R}^n)$. Then we shall say that the distribution $v$ belongs to $H^0_p(\mathbb{R}^n)$, if there exists a function $f \in \mathcal{L}_p(\mathbb{R}^n)$ such that the distributions $v$ and $\mathbf{f}$ coincide, i.e.
$$
\forall \: \varphi \in \mathcal{S}(\mathbb{R}^n) \; \; \; \; v(\varphi) = \int\limits_{\mathbb{R}^n} f \cdot \varphi \; d\mu_n.
$$
The natural norm on the linear space $H^0_p(\mathbb{R}^n)$ is defined as follows:
$$
\forall \: g \in \mathcal{L}_p(\mathbb{R}^n) \quad \| \mathbf{g} \|_{H^0_p(\mathbb{R}^n)} = \| g \|_{\mathcal{L}_p(\mathbb{R}^n)}.
$$
\end{definition}

Actually, Definition \ref{def_H^0_p} allows us to identify the Banach space $\left(H^0_p(\mathbb{R}^n), \: \| \cdot \|_{H^0_p(\mathbb{R}^n)}\right)$ with the Banach space $\left(L_p(\mathbb{R}^n), \: \| \cdot \|_{L_p(\mathbb{R}^n)}\right)$ for arbitrary $p \in [1; \: +\infty)$.

\begin{definition}\label{def_H^s_p}
Let $s \in \mathbb{R}$ and $p \in [1; \: +\infty)$. Let us also introduce the linear operator
$$
J_s \colon \mathcal{S}^{'}(\mathbb{R}^n) \to \mathcal{S}^{'}(\mathbb{R}^n)
$$
as follows: for arbitrary distribution $u \in \mathcal{S}^{'}(\mathbb{R}^n)$ we define
$$
J_s(u) \stackrel{def}{=} \mathcal{F}^{-1}(\varphi_s \cdot \mathcal{F}(u)),
$$
where $\mathcal{F}$ and $\mathcal{F}^{-1}$ are the Fourier transform and the inverse Fourier transform in the dual Schwartz space $\mathcal{S}^{'}(\mathbb{R}^n)$, while the function $\varphi_s \colon \mathbb{R}^n \to \mathbb{C}$ is defined in the following way:
$$
\forall \: x \in \mathbb{R}^n \; \; \: \varphi_s(x) = \left(1 + |x|^2\right)^{\frac{s}{2}}.
$$

Then we define the linear space $H^s_p(\mathbb{R}^n)$ as the set of all distributions $u \in \mathcal{S}^{'}(\mathbb{R}^n)$ such that $J_s(u) \in H^0_p(\mathbb{R}^n)$, with the linear operations being defined naturally. This linear space is equipped with the norm, defined by
$$
\forall \: u \in H^s_p(\mathbb{R}^n) \; \; \; \|u\|_{H^s_p(\mathbb{R}^n)} = \|J_s(u)\|_{H^0_p(\mathbb{R}^n)}.
$$
\end{definition}

Let us fix an arbitrary number $s \in \mathbb{R}$. Then, reasoning by analogy, we can define a linear operator $\mathbb{J}_s \colon \mathcal{S}(\mathbb{R}^n) \to \mathcal{S}(\mathbb{R}^n)$, by defining for an arbitrary function $f \in \mathcal{S}(\mathbb{R}^n)$ the function $J_s(f)$, letting for arbitrary $x \in \mathbb{R}^n$
$$
\left(\mathbb{J}_s(f)\right)(x) \stackrel{def}{=} \left(\mathcal{F}^{-1} \left(\varphi_s \cdot \mathcal{F}(f)\right)\right)(x).
$$
Since not only both the Fourier transform and the inverse Fourier transform but also the operator of multiplication by an arbitrary infinitely differentiable function of at most polynomial growth which never takes $0$ as its value, are linear homeomorphisms of the linear topological space $\mathcal{S}(\mathbb{R}^n)$ onto itself, linear operator $\mathbb{J}_s$ is also a linear homeomorphism of $\mathcal{S}(\mathbb{R}^n)$ onto itself (with respect to $\tau_{\mathcal{S}}$), while linear operator $J_s$ is a linear homeomorphism of $\mathcal{S}'(\mathbb{R}^n)$ onto itself (with respect to $\tau_{\mathcal{S}^{'}}$).

The actions of the operators $J_s \colon \mathcal{S}'(\mathbb{R}^n) \to \mathcal{S}'(\mathbb{R}^n)$ and $\mathbb{J}_s \colon \mathcal{S}(\mathbb{R}^n) \to \mathcal{S}(\mathbb{R}^n)$ are consistent in the following sense. For arbitrary function $f \in \mathcal{S}(\mathbb{R}^n)$ we know that $J_s(\mathbf{f}) \in \mathcal{S}'(\mathbb{R}^n)$ is a regular distribution with the generating function $\mathbb{J}_s(f)$, because we can use the classical equality (which follows immediately from the Plancherel theorem, see, e.g., \cite[Theorem 2.2.14 (5)]{Grafakos_Cl})
$$
\forall \: g \in \mathcal{S}(\mathbb{R}^n), \; \forall \: h \in \mathcal{S}_2(\mathbb{R}^n) \; \; \int\limits_{\mathbb{R}^n} g \cdot h \: d\mu_n = \int\limits_{\mathbb{R}^n} \mathcal{F}(g) \cdot \mathcal{F}^{-1}(h) \: d\mu_n,
$$
and, therefore, we deduce that for any function $\varphi \in \mathcal{S}(\mathbb{R}^n)$ we have
$$
(J_s(\mathbf{f}))(\varphi) = \left(\mathcal{F}^{-1} \left(\varphi_s \cdot \mathcal{F}(\mathbf{f})\right)\right)(\varphi) = \mathbf{f}\left(\mathcal{F}\left(\varphi_s \cdot \mathcal{F}^{-1}(\varphi)\right)\right) =
$$
$$
= \int\limits_{\mathbb{R}^n} f \cdot \mathcal{F}\left(\varphi_s \cdot \mathcal{F}^{-1}(\varphi)\right) \: d\mu_n = \int\limits_{\mathbb{R}^n} \mathcal{F}(f) \cdot \varphi_s \cdot \mathcal{F}^{-1}(\varphi) \: d\mu_n = \int\limits_{\mathbb{R}^n} \mathbb{J}_s(f) \cdot \varphi \: d\mu_n.
$$

Let us also note that, as it is well--known (see, e.g., \cite[Теорема 2.6.1 (a)]{TrBook}), for arbitrary indices $s \in \mathbb{R}$ and $p \in (1; \: +\infty)$ there exists an isometric isomorphism between the spaces
$$
\left(\left(H^s_p(\mathbb{R}^n)\right)^{*}, \: \| \cdot \|_{\left(H^s_p(\mathbb{R}^n)\right)^{*}}\right) \quad \mbox{and} \quad \left(H^{-s}_{p'}(\mathbb{R}^n), \: \| \cdot \|_{H^{-s}_{p'}(\mathbb{R}^n)}\right),
$$
where the former space is defined as a continuous dual space of $\left(H^s_p(\mathbb{R}^n), \: \| \cdot \|_{H^s_p(\mathbb{R}^n)}\right)$.

\begin{definition}\label{loc_unif_spaces}
Let $s \in \mathbb{R}$ and $p \in [1; \: +\infty)$. Let also $\eta \in \mathcal{D}(\mathbb{R}^n)$ be a function which satisfies the conditions
\begin{equation*}\label{eq_eta_1}
a) \; \forall \: x \in \mathbb{R}^n \; \; \; 0 \leqslant \eta(x) \leqslant 1,
\end{equation*}
\begin{equation*}\label{eq_eta_2}
b) \; \forall \: x \in B_1(\mathbf{0}) \; \; \: \eta(x) = 1,
\end{equation*}
\begin{equation*}\label{eq_eta_3}
c) \; \forall \: x \in \mathbb{R}^n \setminus B_2(\mathbf{0}) \; \; \: \eta(x) = 0,
\end{equation*}

Then, firstly, we introduce a linear space
$$
H^s_{p, \: loc}(\mathbb{R}^n) \stackrel{def}{=} \{u \in \mathcal{D}^{'}(\mathbb{R}^n) \: | \; \; \forall \: f \in \mathcal{D}(\mathbb{R}^n) \: \; f \cdot u \in H^s_p(\mathbb{R}^n) \},
$$
and, secondly, we shall say that the distribution $u \in H^s_{p, \: loc}(\mathbb{R}^n)$ is an element of the uniformly localized Bessel potential space $H^s_{p, \: unif, \: \eta}(\mathbb{R}^n)$ if and only if there exists a number $C \in [0; \: +\infty)$ such that for any $z \in \mathbb{R}^n$ we have an estimate
$$
\|\eta_{(z)} \cdot u\|_{H^s_p(\mathbb{R}^n)} \leqslant C.
$$

The set of distributions $H^s_{p, \: unif, \: \eta}(\mathbb{R}^n)$ is a linear subspace of $\mathcal{D}^{'}(\mathbb{R}^n)$, endowed with the norm $\| \cdot \|_{H^s_{p, \: unif, \: \eta}(\mathbb{R}^n)} \colon H^s_{p, \: unif, \: \eta}(\mathbb{R}^n) \to \mathbb{R}$ which is defined as follows:
$$
\forall \: u \in H^s_{p, \: unif, \: \eta}(\mathbb{R}^n) \; \; \; \|u\|_{H^s_{p, \: unif, \: \eta}(\mathbb{R}^n)} \stackrel{def}{=} \sup\limits_{z \in \mathbb{R}^n} \|\eta_{(z)} \cdot u\|_{H^s_p(\mathbb{R}^n)}.
$$
\end{definition}

Since for two arbitrary functions $\eta_1 \in \mathcal{D}(\mathbb{R}^n)$ and $\eta_2 \in \mathcal{D}(\mathbb{R}^n)$, which satisfy the conditions $a)$, $ b)$ and $c)$ from Definition \ref{loc_unif_spaces}, the sets $H^s_{p, \: unif, \: \eta_1}(\mathbb{R}^n)$ and $H^s_{p, \: unif, \: \eta_2}(\mathbb{R}^n)$ do coincide and their natural norms $\| \cdot \|_{s, \: p, \: unif, \eta_1}$ and $\| \cdot \|_{s, \: p, \: unif, \eta_2}$ are equivalent (see an analogous fact in, e.g., \cite[Remark 2.20]{TrBook3}), in the sequel we shall omit the index $\eta$ and use the notation $H^s_{p, \: unif}(\mathbb{R}^n)$ for the spaces of uniformly localized Bessel potential spaces and $\| \cdot \|_{H^s_{p, \: unif}(\mathbb{R}^n)}$ for their norms.

\begin{remark}\label{unif_embeddings}
Using the generalization of the classical Sobolev embedding theorem for the case of the Bessel potential spaces (see, e.g., \cite[Теорема 2.8.1 и замечание 2.8.1.2]{TrBook}) and the fact that a multiplication by an infinitely differentiable function, bounded on $\mathbb{R}^n$ with all its derivatives, defines for $\gamma \in \mathbb{R}$ and $r_1 \geqslant r_2 \geqslant 1$ a linear operator, acting from $H^{\gamma}_{r_1}(\mathbb{R}^n)$ into $H^{\gamma}_{r_2}(\mathbb{R}^n)$ and bounded with respect to the norms $\| \cdot \|_{H^{\gamma}_{r_1}(\mathbb{R}^n)}$ and $\| \cdot \|_{H^{\gamma}_{r_2}(\mathbb{R}^n)}$, it can be shown that for arbitrary $s \in \mathbb{R}, \: t \in \mathbb{R}, \: p \in (1; \: +\infty), \: q \in (1; \: +\infty)$ the following continuous embedding holds true:
$$
\mbox{if} \; \; 1) \; \; p \geqslant q \; \; \mbox{and} \; \; s \geqslant t, \quad \mbox{or} \quad 2) \; \; p \leqslant q \; \; \mbox{and} \; \; s - \frac{n}{p} \geqslant t - \frac{n}{q} \: ,
$$
$$
\mbox{then} \quad H^s_{p, \: unif}(\mathbb{R}^n) \underset{\to}{\subset} H^t_{q, \: unif}(\mathbb{R}^n).
$$
The following well--known (see, e.g., \cite[p. 339]{Marschall}) continuous embedding will also be useful in the sequel:
$$
H^s_p(\mathbb{R}^n) \underset{\to}{\subset} H^s_{p, \: unif}(\mathbb{R}^n).
$$

\end{remark}

\bigskip

\bigskip

\bigskip

\bigskip

{\centerline {\Large 3. Preliminaries: results about multipliers acting between}}

\medskip

{\centerline {\Large the Bessel potential spaces and synopsis }}

\medskip

{\centerline {\Large on the special class of regular distributions.}}

\bigskip

\bigskip

Let us now define the space of multipliers acting in the scale of the Bessel potential spaces.

\begin{definition}\label{def_multipliers}
Let $s \in \mathbb{R}, \: t \in \mathbb{R}, \;  p \in (1; \: +\infty), \: q \in (1; \: +\infty)$. A distribution $\mu \in H^t_{q, \: loc}(\mathbb{R}^n)$ is said to be a multiplier from the space $H^s_p(\mathbb{R}^n)$ to the space $H^t_q(\mathbb{R}^n)$ if there exists a constant $C \in [0; \: +\infty)$ such that for arbitrary function $f \in \mathcal{D}(\mathbb{R}^n)$ the estimate
$$
\|f \cdot \mu\|_{H^t_q(\mathbb{R}^n)} \leqslant C \cdot \|\mathbf{f}\|_{H^s_p(\mathbb{R}^n)}
$$
holds true.

The set $M[H^s_p(\mathbb{R}^n) \to H^t_q(\mathbb{R}^n)]$ of all multipliers from $H^s_p(\mathbb{R}^n)$ to $H^t_q(\mathbb{R}^n)$ is a linear subspace of $H^t_{q, \: loc}(\mathbb{R}^n)$ and it can be endowed with the norm, defined by
$$
\forall \: \mu \in M[H^s_p(\mathbb{R}^n) \to H^t_q(\mathbb{R}^n)]
$$
$$
\|\mu\|_{M[H^s_p(\mathbb{R}^n) \to H^t_q(\mathbb{R}^n)]} = \inf\{C \in [0; \: +\infty) \: | \; \forall \: f \in \mathcal{D}(\mathbb{R}^n) \; \|f \cdot \mu\|_{H^t_q(\mathbb{R}^n)} \leqslant C \cdot \|\mathbf{f}\|_{H^s_p(\mathbb{R}^n)}\}.
$$
\end{definition}

Since for arbitrary $s \in \mathbb{R}$ and $p \in (1; \: +\infty)$ the set of all regular distributions, generated by the functions from $\mathcal{D}(\mathbb{R}^n)$, is dense in $H^s_p(\mathbb{R}^n)$ with respect to the topology induced by the norm $\| \cdot \|_{H^s_p(\mathbb{R}^n)}$ (see a more general result in \cite[Theorem 2.3.2 (b)]{TrBook}), under the conditions of Definition \ref{def_multipliers} a multiplier $u \in M[H^s_p(\mathbb{R}^n) \to H^t_q(\mathbb{R}^n)]$ generates a uniquely defined linear operator $M_u \colon H^s_p(\mathbb{R}^n) \to H^t_q(\mathbb{R}^n)$, which is bounded with respect to the norms $\| \cdot \|_{H^s_p(\mathbb{R}^n)}$ and $\| \cdot \|_{H^t_q(\mathbb{R}^n)}$, with this operator meeting the following conditions:
$$
\forall \: f \in \mathcal{D}(\mathbb{R}^n) \; \; \; M_u(\mathbf{f}) = f \cdot u  \quad \mbox{and} \quad \| u \|_{M[H^s_p(\mathbb{R}^n) \to H^t_q(\mathbb{R}^n)]} = \| M_u \|_{\mathcal{B}(H^s_p(\mathbb{R}^n) \to H^t_q(\mathbb{R}^n))},
$$
where $\| \cdot \|_{\mathcal{B}(H^s_p(\mathbb{R}^n) \to H^t_q(\mathbb{R}^n))}$ is a standard operator norm on the space $\mathcal{B}(H^s_p(\mathbb{R}^n) \to H^t_q(\mathbb{R}^n))$ of all linear operators from $H^s_p(\mathbb{R}^n)$ to $H^t_q(\mathbb{R}^n)$ which are bounded with respect to the standard norms of the respective spaces.

Let us now remind some key facts about the multiplier spaces $M[H^s_p(\mathbb{R}^n) \to H^t_q(\mathbb{R}^n)]$, mostly following \cite{BelShk2017} in the notation and the exposition style.

\bigskip

\begin{proposition}\label{multiplers_symmetricity} $\mathrm{(see \; \cite[Proposition \: 1]{BelShk2017})}$ Let $s \in \mathbb{R}, \: t \in \mathbb{R}, \: p \in (1; \: +\infty)$ and $q \in (1; \: +\infty)$.

\medskip

Then the multiplier spaces $M[H^s_p(\mathbb{R}^n) \to H^t_q(\mathbb{R}^n)]$ and $M[H^{-t}_{q'}(\mathbb{R}^n) \to H^{-s}_{p'}(\mathbb{R}^n)]$ coincide with their norms $\| \cdot \|_{M[H^s_p(\mathbb{R}^n) \to H^t_q(\mathbb{R}^n)]}$ and $\| \cdot \|_{M[H^{-t}_{q'}(\mathbb{R}^n) \to H^{-s}_{p'}(\mathbb{R}^n)]}$ being equivalent.
\end{proposition}

\bigskip

\begin{proposition}\label{easy_multiplers_embedding} $\mathrm{(see \; \cite[Proposition \: 2]{BelShk2017})}$
Let $s \in \mathbb{R}, \: t \in \mathbb{R}, \: p \in (1; \: +\infty)$ and $q \in (1; \: +\infty)$.

Then the following continuous embedding
$$
M[H^s_p(\mathbb{R}^n) \to H^t_q(\mathbb{R}^n)] \underset{\to}{\subset} H^t_{q, \: unif}(\mathbb{R}^n) \cap H^{-s}_{p', \: unif}(\mathbb{R}^n)
$$
is valid, where the natural norm of the space $H^t_{q, \: unif}(\mathbb{R}^n) \cap H^{-s}_{p', \: unif}(\mathbb{R}^n)$ is defined as follows:

for arbitrary $u \in H^t_{q, \: unif}(\mathbb{R}^n) \cap H^{-s}_{p', \: unif}(\mathbb{R}^n)$ we let
$$
\| u \|_{H^t_{q, \: unif}(\mathbb{R}^n) \cap H^{-s}_{p', \: unif}(\mathbb{R}^n)} = \max\left(\| u \|_{H^t_{q, \: unif}(\mathbb{R}^n)}, \; \| u \|_{H^{-s}_{p', \: unif}(\mathbb{R}^n)} \right).
$$
\end{proposition}

\bigskip

\begin{proposition}\label{multiplicative_criterion_for_unif_embedding_into_multipliers} $\mathrm{(An \; immediate \; corollary \; of \; \cite[Lemma \:3, \; Lemma \: 4 \; and \; Lemma \: 5]{BelShk2017})}$
Let $s \in \mathbb{R}$, $\: t \in \mathbb{R}, \: \gamma \in \mathbb{R}, \: r \in (1; \:+\infty)$, $\: p \in (1; \: +\infty)$ and $q \in (1; \: +\infty)$. Let also $p \leqslant q$.

Then the continuous embedding
$$
H^{\gamma}_{r, \: unif}(\mathbb{R}^n) \underset{\to}{\subset} M[H^s_p(\mathbb{R}^n) \to H^t_q(\mathbb{R}^n)]
$$
is valid if and only if there exists a constant $C \in (0; \: +\infty)$, such that for any functions $f \in \mathcal{D}(\mathbb{R}^n)$ and $g \in \mathcal{D}(\mathbb{R}^n)$ the following multiplicative estimate holds true:
$$
\| f \cdot \mathbf{g} \|_{H^t_q(\mathbb{R}^n)} \leqslant C \cdot \| \mathbf{f} \|_{H^{\gamma}_r(\mathbb{R}^n)} \cdot \| \mathbf{g} \|_{H^s_p(\mathbb{R}^n)}.
$$
\end{proposition}

These results are complemented by the description theorems for multiplier spaces, obtained by the author and A.\,A.~Shkalikov in \cite{BelShk2017} and \cite{BelShk2018}.

Firstly, let us formulate the description theorem for multiplier spaces $M[H^s_p(\mathbb{R}^n) \to H^t_q(\mathbb{R}^n)]$ in the case when both smoothness indices $s$ and $t$ are nonnegative and the Strichartz--type conditions are valid.

\begin{theorem}\label{multipliers_description_positive_smoothness} $\mathrm{(\cite[Theorem \: 1]{BelShk2017})}$
Let $p \in (1; \: +\infty), \: q \in (1; \: +\infty), \; s \in \left(\frac{n}{p}; \: +\infty\right), \: t \in [0; \: +\infty)$ and let also $p \leqslant q$ and $s - \frac{n}{p} \geqslant t - \frac{n}{q} \,$.

\medskip

Then we have the set--theoretic equality
$$
M[H^s_p(\mathbb{R}^n) \to H^t_q(\mathbb{R}^n)] = H^t_{q, \: unif}(\mathbb{R}^n)
$$
and the natural norms of these spaces are equivalent.
\end{theorem}

It is worthy to note that the assumptions of Theorem \ref{multipliers_description_positive_smoothness} can be divided into two types: one of these types consists of conditions that guarantee the validity of the continuous embedding $H^s_p(\mathbb{R}^n) \underset{\to}{\subset} H^t_q(\mathbb{R}^n)$, i.e. $p\leqslant q$ and $s - \frac{n}{p} \geqslant t - \frac{n}{q} \,$, and the second of these types is the so--called Strichartz--type condition $s > \frac{n}{p} \,$, importance of which for the problem of finding a constructive description of the multiplier space $M[H^s_p(\mathbb{R}^n) \to H^s_p(\mathbb{R}^n)]$ was outlined by R.\,S.~Strichartz in \cite[Corollary 2.2]{Str67}.

While the problem of investigating the multiplier space $M[H^s_p(\mathbb{R}^n) \to H^t_q(\mathbb{R}^n)]$ in the case when a Strichartz--type condition does not hold will be treated later, it can be easily seen that the conditions of the first type are strictly necessary if we consider the problem of establishing bilateral embeddings of some uniformly localized Bessel potential spaces into the space $M[H^s_p(\mathbb{R}^n) \to H^t_q(\mathbb{R}^n)]$. Indeed, let us assume that for some $\gamma \in \mathbb{R}$ and $r \in (1; \; +\infty)$ there holds a continuous embedding
$$
H^{\gamma}_{r, \: unif}(\mathbb{R}^n) \underset{\to}{\subset} M[H^s_p(\mathbb{R}^n) \to H^t_q(\mathbb{R}^n)].
$$
Then we can consider the regular distribution $\mathbf{Id}$, generated by the identical function $Id_{\mathbb{R}^n}$, which belongs to $\mathcal{L}_{1, \: loc}(\mathbb{R}^n)$ and has at most polynomial growth on $\mathbb{R}^n$, and, since it can be easily shown that $\mathbf{Id} \in H^{\gamma}_{r, \: unif}(\mathbb{R}^n)$, then for arbitrary $u \in H^s_p(\mathbb{R}^n)$ we have $M_{\mathbf{Id}}(u) = u \in H^t_q(\mathbb{R}^n)$ and for some constant $C \in (0; \; +\infty)$, not depending on the choice of $u \in H^s_p(\mathbb{R}^n)$, we have the estimate
$$
\| u \|_{H^t_q(\mathbb{R}^n)} \leqslant C \cdot \| u \|_{H^s_p(\mathbb{R}^n)}.
$$
This means exactly that the continuous embedding $H^s_p(\mathbb{R}^n) \underset{\to}{\subset} H^t_q(\mathbb{R}^n)$ holds true, which implies that the conditions $p \leqslant q$ and $s - \frac{n}{p} \geqslant t - \frac{n}{q}$ are valid.

The description theorem for multiplier spaces $M[H^s_p(\mathbb{R}^n) \to H^t_q(\mathbb{R}^n)]$ in the case when smoothness indices $s$ and $t$ have different signs but the Strichartz--type conditions are still valid states the following.

\begin{theorem}\label{multipliers_description_different_sign_smoothness_indices} $\mathrm{(\cite[Theorem \: 1]{BelShk2018})}$
Let $p \in (1; \: +\infty), \: q \in (1; \: +\infty), \: p \leqslant q, \; s \in [0; \: +\infty), \: t \in [0; \: +\infty)$ and let one of the following conditions be valid:

either $1) \; s \geqslant t, \: s > \frac{n}{p}\,$;

or $2) \; t \geqslant s, \: t > \frac{n}{q'}\,$.

Then we have the set--theoretic equality
$$
M[H^s_p(\mathbb{R}^n) \to H^{-t}_q(\mathbb{R}^n)] = H^{-t}_{q, \: unif}(\mathbb{R}^n) \cap H^{-s}_{p', \: unif}(\mathbb{R}^n)
$$
and the natural norms of these spaces (with the norm on $H^{-t}_{q, \: unif}(\mathbb{R}^n) \cap H^{-s}_{p', \: unif}(\mathbb{R}^n)$ being defined as a maximum of the norms $\| \cdot \|_{H^s_p(\mathbb{R}^n)}$ and $\| \cdot \|_{H^t_q(\mathbb{R}^n)}$) are equivalent.
\end{theorem}

Once again, just like this was the case with the assumptions of Theorem \ref{multipliers_description_positive_smoothness}, the condition $p \leqslant q$ from the assumptions of Theorem \ref{multipliers_description_different_sign_smoothness_indices} is necessary for the continuous embedding $H^{\gamma}_{r, \: unif}(\mathbb{R}^n) \underset{\to}{\subset} M[H^s_p(\mathbb{R}^n) \to H^{-t}_q(\mathbb{R}^n)]$ to hold true for some $\gamma \in \mathbb{R}$ and $r \in (1; \; +\infty)$ (we can also see that the condition $s - \frac{n}{p} \geqslant -t - \frac{n}{q}$ holds true in the situation when $s > \frac{n}{p}$ as $-t - \frac{n}{q} < 0$).

So, in the situation when the Strichartz--type conditions on the smoothness indices do not hold, we are interested in finding bilateral embeddings of the type
$$
H^{\gamma_1}_{r_1, \: unif}(\mathbb{R}^n) \underset{\to}{\subset} M[H^s_p(\mathbb{R}^n) \to H^t_q(\mathbb{R}^n)] \underset{\to}{\subset} H^{\gamma_2}_{r_2, \: unif}(\mathbb{R}^n)
$$
and proving the optimal character of the lower indices in the model situation $p = q = 2$, which allows us to understand the significance of the Strichartz--type conditions for the constructive description of the multiplier spaces $M[H^s_p(\mathbb{R}^n) \to H^t_q(\mathbb{R}^n)]$ in terms of the scale of the uniformly localized Bessel potential spaces.

These bilateral embeddings for the situation when the smoothness indices are of different signs were obtained in \cite{NZSh2006}. A slight reformulation of \cite[Lemma 3]{NZSh2006} and \cite[Lemma 6]{NZSh2006} leads us to the following statement.

\begin{proposition}\label{different_sign_smoothness_indices_bilateral_embeddings}
Let $s \in [0; \: +\infty], \; t \in [0; \: +\infty)$ and $\max(s; \: t) \in \left(0; \: \frac{n}{2}\right)$.

Then the following continuous embeddings hold true:
$$
H^{-\min(s; \: t)}_{\frac{n}{\max(s; \: t)}, \: unif}(\mathbb{R}^n) \underset{\to}{\subset} M[H^s_2(\mathbb{R}^n) \to H^{-t}_2(\mathbb{R}^n)]
$$
\centerline{and}
$$
M[H^s_2(\mathbb{R}^n) \to H^{-t}_2(\mathbb{R}^n)] \underset{\to}{\subset} H^{-\min(s; \: t)}_{2, \: unif}(\mathbb{R}^n).
$$
\end{proposition}

In the paper \cite{BelShkArXiv2019} the exact character of the left embedding in Proposition \ref{different_sign_smoothness_indices_bilateral_embeddings} was established by considering the regular distributions $\left(\mathbf{f}_{\alpha} \: | \; \alpha \in (0; \: n)\right)$, where for any $\alpha \in (0; \: n)$ the regular distribution $\mathbf{f}_{\alpha}$ is generated by the function $f_{\alpha}$, with this function being defined as follows:
$$
\forall \: x \in \mathbb{R}^n \; \; \: f_{\alpha}(x) = \begin{cases}
\frac{1}{\left|x\right|^{\alpha}}, \; \: \mbox{if} \: \: x \neq \mathbf{0};\\
0, \; \: \mbox{if} \: \: x = \mathbf{0}.
\end{cases}
$$

The main result of that paper is the following one.

\begin{theorem}\label{exact_embedding} $\mathrm{\cite[Theorem \: 3]{BelShkArXiv2019}}$
Let $s \in [0; \: +\infty), \: t \in [0; \: +\infty)$ and $\max (s, t) \in \left(0; \: \frac{n}{2}\right)$. Then for arbitrary number $\varepsilon \in \left(0; \: \frac{n}{\max (s, t)} - 2\right)$, there exists a positive number
$$
\delta(\varepsilon) \in \Bigl(0; \frac{n}{2} - \max (s, t)\Bigr),
$$
such that for $\alpha(\varepsilon) = s + t + \delta(\varepsilon)$ we have
$$
\mathbf{f_{\alpha(\varepsilon)}} \in H^{-\min (s, t)}_{\frac{n}{\max (s, t)} - \varepsilon, \: unif}(\mathbb{R}^n) \setminus \: M[H^s_2(\mathbb{R}^n) \to H^{-t}_2(\mathbb{R}^n)].
$$
\end{theorem}

Let us also briefly recall when this regular distribution $\mathbf{f}_{\alpha}$ can be seen as a well-defined element of the dual Schwartz space $\mathcal{S}^{'}(\mathbb{R}^n)$, with our presentation mostly following \cite{BelShkArXiv2019}.

First of all, let us note that since for any numbers $r_1 \in (0; \: +\infty)$ and $r_2 \in (0; \: +\infty)$, such that $r_1 > r_2$, we have
$$
\int\limits_{B_{r_1, \, r_2}(\mathbf{0})} f_{\alpha}(x) \: d\mu_n(x) = C(n) \cdot \int\limits_{[r_1, \, r_2]} \frac{1}{r^{\alpha}} \cdot r^{n - 1} \: dr = C(n) \cdot \int\limits_{[r_1, \, r_2]} \frac{1}{r^{\alpha - n + 1}} \: dr,
$$
where $C(n)$ is a constant equal to the surface area of the $(n-1)$--dimensional hypersphere, i.e.
$$
C(n) = \frac{2 \cdot \pi^{\frac{n}{2}}}{\Gamma(\frac{n}{2})} \: .
$$
Therefore, as for  $\alpha < n$ we have $f_{\alpha} \in \mathcal{L}_{1, \: loc}(\mathbb{R}^n)$ and, since $f_{\alpha}$ has at most polynomial growth outside $B_1(\mathbf{0})$, in this case we can define a regular functional $\mathbf{f}_{\alpha}$ on $\mathcal{S}(\mathbb{R}^n)$ and, moreover, $\mathbf{f}_{\alpha} \in \mathcal{S}'(\mathbb{R}^n)$.

On the other hand, if $\alpha \geqslant n$, then the function $f_{\alpha}$ is not integrable on $B_1(\mathbf{0})$ and, taking an arbitrary function $\varphi \in \mathcal{D}(\mathbb{R}^n)$ that is equal to $1$ on $B_1(\mathbf{0})$, we arrive at the fact that the value $\mathbf{f}_{\alpha}(\varphi)$ can not be correctly defined.

Thereby, for $\alpha \in (0; \: +\infty)$ the condition $\alpha < n$ is a criterion for the regular distribution $\mathbf{f}_{\alpha} \colon \mathcal{S}(\mathbb{R}^n) \to \mathbb{C}$ to be well--defined on $\mathcal{S}(\mathbb{R}^n)$ and, moreover, for $\alpha \in (0; \: n)$ this distribution is an element of the dual Schwartz space $\mathcal{S}^{'}(\mathbb{R}^n)$.

Also in the sequel we shall need the following fact, which constitutes a partial case of a classical fact from the theory of Sobolev--type spaces (see, e.g., \cite[Lemma 2.3.1]{RSbook}).

\begin{proposition}\label{Sickel_Runst_special_distribution}
Let $p \in [1; \: +\infty), \: s > 0, \: \alpha \in (0; \: n)$ and let also a function $\eta \in \mathcal{D}(\mathbb{R}^n)$ satisfy the following conditions:
\begin{equation*}\label{eq_eta_1_new}
a) \; \forall \: x \in \mathbb{R}^n \; \; \; 0 \leqslant \eta(x) \leqslant 1,
\end{equation*}
\begin{equation*}\label{eq_eta_2_new}
b) \; \forall \: x \in B_1(\mathbf{0}) \; \; \: \eta(x) = 1,
\end{equation*}
\begin{equation*}\label{eq_eta_3_new}
c) \; \forall \: x \in \mathbb{R}^n \setminus B_2(\mathbf{0}) \; \; \: \eta(x) = 0.
\end{equation*}

Then
$$
\eta \cdot \mathbf{f_{\alpha}} \in H^s_p(\mathbb{R}^n)\quad \mbox{if and only if} \quad \alpha < \frac{n}{p} - s.
$$
\end{proposition}

\bigskip

\bigskip

\bigskip

\bigskip

{\centerline {\Large 4. Multipliers between two Bessel potential spaces with the}}
{\centerline {\Large positive smoothness indices: bilateral continuous}}
{\centerline {\Large  embeddings and their exact character.}}

\bigskip

\bigskip

Let us consider the multiplier space $M[H^s_2(\mathbb{R}^n) \to H^t_2(\mathbb{R}^n)]$ in the situation when
$$
s \in [0; \: +\infty), \; t \in [0; \: +\infty) \; \: \mbox{and} \; \: \max(s; \: t) \in \left(0; \: \frac{n}{2}\right).
$$
Then the first type assumptions from Theorem \ref{multipliers_description_positive_smoothness} turn into a tautological inequality $2 \leqslant 2$ and an inequality $s - \frac{n}{2} \geqslant t - \frac{n}{2}$, which implies $s \geqslant t$. This assumption is not really restrictive since it is well--known  that if $s < t$, then the multiplier space $M[H^s_2(\mathbb{R}^n) \to H^t_2(\mathbb{R}^n)]$ consists only of a trivial multiplier, which maps any element of $H^s_p(\mathbb{R}^n)$ into the identically zero distribution.

We shall also need the following multiplicative estimate in the scale of the Lizorkin--Triebel spaces.

\begin{proposition}\label{Sickel_Runst_multiplicative_estimate} $\mathrm{(see \; a \; more \; general \; fact \; in \: \cite[Theorem \: 4.4.2]{RSbook})}$
Let
$$
\widetilde{s_1} \in (0; \: +\infty), \: \widetilde{s_2} \in (0; \: +\infty), \: \widetilde{s_2} \geqslant \widetilde{s_1}, \: \widetilde{p} \in [1; \: +\infty), \: \widetilde{p_1} \in [1; \; +\infty), \: \widetilde{p_2} \in [1; \: +\infty),
$$
$$
\widetilde{q} \in (1; \: +\infty), \: \widetilde{q_1} \in (1; \; +\infty) \; \; \mbox{and} \; \; \widetilde{q_2} \in (1; \: +\infty).
$$
Let us also assume that the following conditions hold true
$$
1) \; \frac{1}{\widetilde{p}} \leqslant \frac{1}{\widetilde{p_1}} + \frac{1}{\widetilde{p_2}}\, ;
$$
$$
2) \; \widetilde{s_1} + \widetilde{s_2} > \frac{n}{\widetilde{p_1}} + \frac{n}{\widetilde{p_2}} - n;
$$
$$
3) \; \widetilde{s_2} > \widetilde{s_1} \; \: \mbox{and} \; \; \widetilde{q} \geqslant \widetilde{q_1}; 
$$
$$
4) \; \max\left(\frac{n}{\widetilde{p_1}} - \widetilde{s_1}; \: \frac{n}{\widetilde{p_2}} - \widetilde{s_2}\right) > 0 \; \; \mbox{and} \; \; \frac{n}{\widetilde{p}} - \widetilde{s_1} = \max\left(\frac{n}{\widetilde{p_1}} - \widetilde{s_1}; \: 0\right) + \max\left(\frac{n}{\widetilde{p_2}} - \widetilde{s_2}; \: 0\right);
$$
$$
5) \; \left\{i \in \{1; \: 2\} \: | \; \widetilde{s_i} = \frac{n}{\widetilde{p_i}} \; \mbox{and} \; \widetilde{p_i} > 1\right\} = \varnothing.
$$
Then there exists a constant $C \in [0; \: +\infty)$, such that for arbitrary functions $f \in \mathcal{D}(\mathbb{R}^n)$ and $g \in \mathcal{D}(\mathbb{R}^n)$ we have the following multiplicative estimate:
$$
\| f \cdot \mathbf{g} \|_{F^{\widetilde{s_1}}_{\widetilde{p}, \: \widetilde{q}}(\mathbb{R}^n)} \leqslant C \cdot \| \mathbf{f} \|_{F^{\widetilde{s_1}}_{\widetilde{p_1}, \: \widetilde{q_1}}(\mathbb{R}^n)} \cdot \| \mathbf{g} \|_{F^{\widetilde{s_2}}_{\widetilde{p_2}, \: \widetilde{q_2}}(\mathbb{R}^n)}.
$$
\end{proposition}

\bigskip

Now we are ready to prove a result which in the non--Strichartz case establishes bilateral continuous embeddings, featuring multiplier space between two Bessel potential spaces with nonnegative smoothness indices and the spaces from the scale of the uniformly localized Bessel potential spaces. This result can be seen as a counterpart to Proposition \ref{different_sign_smoothness_indices_bilateral_embeddings} which covers the case when the smoothness indices have different signs.

\medskip

\begin{theorem}\label{bilateral_embeddings_positive_smoothness_Bessel_potential_spaces}
Let $s \in (0; \: +\infty), \: t \in [0; \: +\infty), \: s > t$ and $s < \frac{n}{2}$. Then the continuous embeddings
$$
H^t_{\frac{n}{s}, \: unif}(\mathbb{R}^n) \underset{\to}{\subset} M[H^s_2(\mathbb{R}^n) \to H^t_2(\mathbb{R}^n)] \quad \mbox{and} \quad M[H^s_2(\mathbb{R}^n) \to H^t_2(\mathbb{R}^n)] \underset{\to}{\subset} H^t_{2, \: unif}(\mathbb{R}^n)
$$
hold true.
\end{theorem}

Proof. By Proposition \ref{easy_multiplers_embedding} we have the continuous embedding
$$
M[H^s_2(\mathbb{R}^n) \to H^t_2(\mathbb{R}^n)] \underset{\to}{\subset} H^t_2(\mathbb{R}^n) \cap H^{-s}_2(\mathbb{R}^n).
$$
Since $t > -s$, by Remark \ref{unif_embeddings} we have the validity of the continuous embedding
$$
H^t_2(\mathbb{R}^n) \underset{\to}{\subset} H^{-s}_2(\mathbb{R}^n)
$$
and, hence, the space $H^t_2(\mathbb{R}^n) \cap H^{-s}_2(\mathbb{R}^n)$ coincides with the space $H^t_2(\mathbb{R}^n)$ with their natural norms being equivalent.

Summing it up, we arrive at the validity of the continuous embedding
$$
M[H^s_2(\mathbb{R}^n) \to H^t_2(\mathbb{R}^n)] \underset{\to}{\subset} H^t_2(\mathbb{R}^n).
$$

Let us prove the other embedding
$$
H^t_{\frac{n}{s}, \: unif}(\mathbb{R}^n) \underset{\to}{\subset} M[H^s_2(\mathbb{R}^n) \to H^t_2(\mathbb{R}^n)].
$$

Firstly, let us consider the situation when $t > 0$.

Since it is a well--known fact (see, e.g., \cite[Theorem 2.3.3]{TrBook}) that for any indices $s \in \mathbb{R}$ and $p \in [1; \: +\infty)$ we have the set--theoretic equality $F^s_{p, \: 2}(\mathbb{R}^n) = H^s_p(\mathbb{R}^n)$ with the norms $\| \cdot \|_{F^s_{p, \: 2}(\mathbb{R}^n)}$ and $\| \cdot \|_{H^s_p(\mathbb{R}^n)}$ being equivalent, we employ Proposition \ref{Sickel_Runst_multiplicative_estimate}, where we take $\widetilde{s_1} = t, \; \widetilde{p} = 2, \; \widetilde{q} = 2, \; \widetilde{p_1} = \frac{n}{s} \, , \; \widetilde{q_1} = 2, \; \widetilde{s_2} = s, \; \widetilde{p_2} = 2, \; \widetilde{q_2} = 2$, we obtain the fact that there exists a constant $C \in [0; \: +\infty)$, such that for arbitrary functions $f \in \mathcal{D}(\mathbb{R}^n)$ and $g \in \mathcal{D}(\mathbb{R}^n)$ we have the estimate
\begin{equation}\label{RuSi_estimate}
\| f \cdot \mathbf{g} \|_{H^t_2(\mathbb{R}^n)} \leqslant C \cdot \| \mathbf{f} \|_{H^s_2(\mathbb{R}^n)} \cdot \| \mathbf{g} \|_{H^t_{\frac{n}{s}}(\mathbb{R}^n)},
\end{equation}
if the following conditions are valid:
$$
1) \: \frac{1}{\widetilde{p}} = \frac{1}{2} \leqslant \frac{1}{\widetilde{p_1}} + \frac{1}{\widetilde{p_2}} =  \frac{s}{n} + \frac{1}{2}; \quad 2) \: \widetilde{s_1} + \widetilde{s_2} = t + s > \frac{n}{\widetilde{p_1}} + \frac{n}{\widetilde{p_2}} - n = s + \frac{n}{2} - n = s - \frac{n}{2};
$$
$$
3) \; \widetilde{s_2} = s > \widetilde{s_1} = t \; \: \mbox{and} \; \: \widetilde{q} = 2 \geqslant \widetilde{q_1} = 2;
$$
$$
4) \; \max\left(\frac{n}{\widetilde{p_1}} - \widetilde{s_1}; \: \frac{n}{\widetilde{p_2}} - \widetilde{s_2}\right) = \max\left(s - t; \frac{n}{2} - s\right) > 0 \; \: \mbox{and}
$$
$$
\frac{n}{\widetilde{p}} - \widetilde{s_1} = \frac{n}{2} - t = s - t + \frac{n}{2} - s = 
$$
$$
= \max(s - t; \: 0) + \max\left(\frac{n}{2} - s; \: 0\right) = \max\left(\frac{n}{\widetilde{p_1}} - \widetilde{s_1} \, ; \: 0\right) + \max\left(\frac{n}{\widetilde{p_2}} - \widetilde{s_2} \, ; \: 0\right);
$$
$$
5) \; \widetilde{s_1} = t \neq \frac{\widetilde{n}}{\widetilde{p_1}} = s \; \; \mbox{and} \; \; \: \widetilde{s_2} = s \neq \frac{\widetilde{n}}{\widetilde{p_2}} = \frac{n}{2}.
$$

Now, we remark that the condition $1)$ holds true since $s > 0$, the condition $2)$ holds true since $t > 0 > -\frac{n}{2}$, the condition $3)$ holds true since $s > t$, the condition $4)$ holds true since $\frac{n}{2} - s > 0$ and, finally, the condition $5)$ holds true since $s > t$ and $s < \frac{n}{2}\,$.

So, we conclude that the estimate \eqref{RuSi_estimate} is valid.

By Proposition \ref{multiplicative_criterion_for_unif_embedding_into_multipliers} (in our situation both $p$ and $q$ are equal to $2$), the validity of the estimate \eqref{RuSi_estimate} guarantees the validity of the continuous embedding
$$
H^t_{\frac{n}{s}, \: unif}(\mathbb{R}^n) \underset{\to}{\subset} M[H^s_2(\mathbb{R}^n) \to H^t_2(\mathbb{R}^n)].
$$

Now only the limit case $t = 0$ is left to consider. In that case the validity of the continuous embedding
$$
H^0_{\frac{n}{s}, \: unif}(\mathbb{R}^n) \underset{\to}{\subset} M[H^s_2(\mathbb{R}^n) \to H^0_2(\mathbb{R}^n)]
$$
follows immediately from Proposition \ref{bilateral_embeddings_positive_smoothness_Bessel_potential_spaces} (see \cite[Lemma 6]{NZSh2006}).

This concludes the proof of Theorem \ref{bilateral_embeddings_positive_smoothness_Bessel_potential_spaces}.

\bigskip

\bigskip

Now we turn to the problem of establishing the exact character of the index $\frac{n}{s}$ in the continuous embedding
$$
H^t_{\frac{n}{s}, \: unif}(\mathbb{R}^n) \underset{\to}{\subset} M[H^s_2(\mathbb{R}^n) \to H^t_2(\mathbb{R}^n)]
$$
for the situation when $0 < t < s < \frac{n}{2}\,$.

\bigskip

\begin{proposition}\label{smoothed_special_functional_sufficient_condition}
Let $t \in \left(-\infty; \: \frac{n}{2}\right)$ and $q \in (1; \: +\infty)$. Let also a function $\eta \in \mathcal{D}(\mathbb{R}^n)$ satisfy the following conditions:
\begin{equation*}\label{eq_eta_1_new}
a) \; \forall \: x \in \mathbb{R}^n \; \; \; 0 \leqslant \eta(x) \leqslant 1,
\end{equation*}
\begin{equation*}\label{eq_eta_2_new}
b) \; \forall \: x \in B_1(\mathbf{0}) \; \; \; \eta(x) = 1,
\end{equation*}
\begin{equation*}\label{eq_eta_3_new}
c) \; \forall \: x \in \mathbb{R}^n \setminus B_2(\mathbf{0}) \; \; \; \eta(x) = 0,
\end{equation*}
Finally, let $\alpha \in \left(0; \: \min\left(n; \: \frac{n}{q} - t\right)\right)$.

\medskip

Then
$$
\eta \cdot \mathbf{f}_{\alpha} \in H^t_{q, \: unif}(\mathbb{R}^n).
$$
\end{proposition}

Proof. Since $\alpha \in (0; \: n)$ it follows that $\mathbf{f}_{\alpha}$ is well--defined as a distribution from $\mathcal{S}^{'}(\mathbb{R}^n)$. Because the multiplication by a function from $\mathcal{D}(\mathbb{R}^n)$ induces a continuous operator in the dual Schwartz space $\mathcal{S}^{'}(\mathbb{R}^n)$, we know that the distribution $\eta \cdot \mathbf{f}_{\alpha}$ is also well--defined as an element of $\mathcal{S}^{'}(\mathbb{R}^n)$.

Then, by Proposition \ref{Sickel_Runst_special_distribution}, we have $\eta \cdot \mathbf{f}_{\alpha} \in H^t_q(\mathbb{R}^n)$. By Remark \ref{unif_embeddings}, the continuous embedding 
$$
H^t_q(\mathbb{R}^n) \underset{\to}\subset H^t_{q, \: unif}(\mathbb{R}^n)
$$
holds true and, therefore, $\eta \cdot \mathbf{f}_{\alpha} \in H^t_{q, \: unif}(\mathbb{R}^n)$.

This concludes the proof of Proposition \ref{smoothed_special_functional_sufficient_condition}.

\bigskip

\begin{proposition}\label{multiplier_necessity_smooth_special_function}
Let $s \in \left(0; \: \frac{n}{2}\right), \: t \in \left(0; \: \frac{n}{2}\right)$ and $s > t$. Let also a function $\eta \in \mathcal{D}(\mathbb{R}^n)$ satisfy the following conditions:
\begin{equation*}\label{eq_eta_1_new}
a) \; \forall \: x \in \mathbb{R}^n \; \; \; 0 \leqslant \eta(x) \leqslant 1,
\end{equation*}
\begin{equation*}\label{eq_eta_2_new}
b) \; \forall \: x \in B_1(\mathbf{0}) \; \; \; \eta(x) = 1,
\end{equation*}
\begin{equation*}\label{eq_eta_3_new}
c) \; \forall \: x \in \mathbb{R}^n \setminus B_2(\mathbf{0}) \; \; \; \eta(x) = 0,
\end{equation*}
Finally, let $\alpha \in \left(0; \: \frac{n}{2}\right)$ and $\eta \cdot \mathbf{f}_{\alpha} \in M[H^s_2(\mathbb{R}^n) \to H^t_2(\mathbb{R}^n)]$.

\medskip

Then $\alpha \leqslant s - t$.
\end{proposition}

\medskip

Proof. Fix an arbitrary $\varepsilon \in \left(0; \: \frac{n}{2} - s\right)$. 

Let
$$
\beta_{\varepsilon} \stackrel{def}{=} \frac{n}{2} - s - \varepsilon.
$$
Then $\beta_{\varepsilon} > 0$ and also $\beta_{\varepsilon} < \frac{n}{2} < n$. So, as $\beta_{\varepsilon} \in (0; \: n)$, the functional $\eta \cdot \mathbf{f}_{\beta_{\varepsilon}}$ is well--defined as an element of $\mathcal{S}^{'}(\mathbb{R}^n)$. Since $\beta_{\varepsilon} < \frac{n}{2} - s$, by Proposition \ref{Sickel_Runst_special_distribution} we have 
$$
\eta \cdot \mathbf{f}_{\beta_{\varepsilon}} \in H^s_2(\mathbb{R}^n).
$$

As $\eta \cdot \mathbf{f}_{\alpha} \in M[H^s_2(\mathbb{R}^n) \to H^t_2(\mathbb{R}^n)]$, it follows that
$$
M_{\eta \cdot \mathbf{f}_{\alpha}}(\eta \cdot \mathbf{f}_{\beta_{\varepsilon}}) \in H^t_2(\mathbb{R}^n).
$$

Because $f_{\alpha} \cdot f_{\beta_{\varepsilon}} = f_{\alpha + \beta_{\varepsilon}}$ and $\alpha + \beta_{\varepsilon} < \frac{n}{2} + \frac{n}{2} = n$, it follows 
that 
$$
M_{\eta \cdot \mathbf{f}_{\alpha}}(\eta \cdot \mathbf{f}_{\beta_{\varepsilon}}) = \eta^2 \cdot \mathbf{f}_{\alpha + \beta_{\varepsilon}}.
$$
As the function $\eta^2 \in \mathcal{D}(\mathbb{R}^n)$ satisfies the same conditions $a), \: b)$ and $c)$ as $\eta$, once again using Proposition \ref{Sickel_Runst_special_distribution}, we obtain that the condition $\eta^2 \cdot \mathbf{f}_{\alpha + \beta_{\varepsilon}} \in H^t_2(\mathbb{R}^n)$ implies the validity of the estimate $\alpha + \beta_{\varepsilon} < \frac{n}{2} - t$.

Therefore,
\begin{equation}\label{technical_eq_2}
\alpha < \frac{n}{2} - t - \beta_{\varepsilon} = s - t + \varepsilon.
\end{equation}

Since $\varepsilon$ was taken arbitrarily from the interval $\left(0; \: \frac{n}{2} - s\right)$, taking the limits in \eqref{technical_eq_2} as $\varepsilon \to 0_+$, we obtain the needed estimate $\alpha \leqslant s - t$.

This concludes the proof of Proposition \ref{multiplier_necessity_smooth_special_function}.

\bigskip

\bigskip

\begin{theorem}\label{positive_smoothness_exact_embedding}
Let $s \in \left(0; \: \frac{n}{2}\right), \: t \in \left(0; \: \frac{n}{2}\right)$ and $s > t$. Then for any number $\varepsilon \in \left(0; \; \frac{n}{s} - 2\right)$ there exists a distribution
$$
u_{\varepsilon} \in H^t_{\frac{n}{s} - \varepsilon, \: unif}(\mathbb{R}^n) \setminus M[H^s_2(\mathbb{R}^n) \to H^t_2(\mathbb{R}^n)].
$$
\end{theorem}

Proof. Fix an arbitrary number $\varepsilon_0 \in \left(0; \; \frac{n}{s} - 2\right)$.

Also fix an arbitrary function $\eta \in \mathcal{D}(\mathbb{R}^n)$, such that
\begin{equation*}\label{eq_eta_1_new}
a) \; \forall \: x \in \mathbb{R}^n \; \; \; 0 \leqslant \eta(x) \leqslant 1,
\end{equation*}
\begin{equation*}\label{eq_eta_2_new}
b) \; \forall \: x \in B_1(\mathbf{0}) \; \; \; \eta(x) = 1,
\end{equation*}
\begin{equation*}\label{eq_eta_3_new}
c) \; \forall \: x \in \mathbb{R}^n \setminus B_2(\mathbf{0}) \; \; \; \eta(x) = 0.
\end{equation*}

Let
$$
\delta(\varepsilon_0) \stackrel{def}{=} \frac{\varepsilon_0 \cdot s^2}{2 \cdot (n - \varepsilon_0 \cdot s)} \quad \mbox{and} \quad \alpha(\varepsilon_0) \stackrel{def}{=} s - t + \delta(\varepsilon_0).
$$
Since $\varepsilon_0 > 0, \ s > 0$ and $n - \varepsilon_0 \cdot s > n - n + 2 \cdot s = 2 \cdot s > 0$, it follows that $\delta(\varepsilon_0) > 0$ and, consequently, $\alpha(\varepsilon_0) > s - t$.

On the other hand,
$$
\alpha(\varepsilon_0) = s - t + \delta(\varepsilon_0) = s - t + \frac{\varepsilon_0 \cdot s^2}{2 \cdot (n - \varepsilon_0 \cdot s)} <
$$
$$
< s - t + \frac{\varepsilon_0 \cdot s}{\frac{n}{s} - \varepsilon_0} = \frac{n}{\frac{n}{s} - \varepsilon_0} - t < \frac{n}{\frac{n}{s} - \varepsilon_0} < \frac{n}{2} \: .
$$

Therefore, by Proposition \ref{multiplier_necessity_smooth_special_function}, it follows that $\eta \cdot \mathbf{f}_{\alpha(\varepsilon_0)} \notin M[H^s_2(\mathbb{R}^n) \to H^t_2(\mathbb{R}^n)]$.

Then, having in mind the previous chain of inequalities, we know that
$$
\alpha(\varepsilon_0) < \frac{n}{\frac{n}{s} - \varepsilon_0} - t,
$$
which, combined with the estimate $\alpha(\varepsilon_0) < \frac{n}{2} < n$, allows us to employ Proposition \ref{smoothed_special_functional_sufficient_condition} which implies
$$
\eta \cdot \mathbf{f}_{\alpha(\varepsilon_0)} \cdot \in H^t_{\frac{n}{s} - \varepsilon_0, \: unif}(\mathbb{R}^n).
$$

So, taking our regular distribution $\eta \cdot \mathbf{f}_{\alpha(\varepsilon_0)}$ as $u_{\varepsilon_0}$, we obtain that $u_{\varepsilon_0}$ belongs to the uniformly localized Bessel potential space $H^t_{\frac{n}{s} - \varepsilon_0, \: unif}(\mathbb{R}^n)$ but not to the multiplier space $M[H^s_2(\mathbb{R}^n) \to H^t_2(\mathbb{R}^n)]$. 

Since $\varepsilon_0 \in \left(0; \; \frac{n}{s} - 2\right)$ was taken arbitrarily, this gives us a required result.

This concludes the proof of Theorem \ref{positive_smoothness_exact_embedding}.

\bigskip

Addresses:

A.\,A.~Belyaev,
Lomonosov Moscow State University,
Department of Mechanics and Mathematics;
Peoples' Friendship University of Russia,
S.\,M.~Nikol'skii Mathematical Institute
{email: belyaev\_aa@pfur.ru}

\end{document}